\documentclass[11pt,twoside]{article}
\usepackage{amsmath,amsfonts,amssymb,theorem}

\newtheorem{thm}{Theorem}[section] {
  
  \newtheorem{lemm}[thm]{Lemma}
  
  \newtheorem{cor}[thm]{Corollary}
  \newtheorem{defi}[thm]{Definition}
{\theorembodyfont{\rmfamily}  \newtheorem{rem}[thm]{Remark} }
 }

\def\Xint#1{\mathchoice
   {\XXint\displaystyle\textstyle{#1}}%
   {\XXint\textstyle\scriptstyle{#1}}%
   {\XXint\scriptstyle\scriptscriptstyle{#1}}%
   {\XXint\scriptscriptstyle\scriptscriptstyle{#1}}%
   \!\int}
\def\XXint#1#2#3{{\setbox0=\hbox{$#1{#2#3}{\int}$}
     \vcenter{\hbox{$#2#3$}}\kern-.5\wd0}}

\def\aint{\Xint\diagup}

\newcommand{\bke}[1]{\left( #1 \right)}
\newcommand{\bkt}[1]{\left[ #1 \right]}
\newcommand{\bket}[1]{\left\{ #1 \right\}}
\newcommand{\norm}[1]{\left\Vert #1 \right\Vert}
\newcommand{\abs}[1]{\left| #1 \right|}

\newcommand{\la}{\lambda}
\newcommand{\e}{\epsilon}
\newcommand{\eps}{{\epsilon}}
\renewcommand{\div}{\mathop{\mathrm{div}}}
\newcommand{\curl}{\mathop{\mathrm{curl}}}
\newcommand{\si}{\sigma}
\renewcommand{\th}{\theta}
\newcommand{\R}{{\mathbb R }}
\newcommand{\calS}{{\cal S }}
\newcommand{\calC}{{\mathcal C}}
\newcommand{\calP}{{\mathcal P}}
\newcommand{\pr}{{\partial}}
\newcommand{\Om}{{\Omega}}
\newcommand{\nb}{{\nabla }}
\newcommand{\ind}{{\,\mathrm{d}}}
\newcommand{\qed}{\hfill\fbox{}\par\vspace{0.3mm}}
\newcommand{\ga}{{\gamma}}
\newcommand{\MAINPF}[1]{\par \medskip
    \noindent{\bf Proof of Theorem \ref{TH1} #1.}}
\newenvironment{pf}{{\bf Proof.}} {\hfill\qed}

\begin{document}

\title{Interior regularity criteria for suitable weak solutions
of the Navier-Stokes equations}

\author{Stephen Gustafson$^1$, Kyungkeun Kang$^2$, and Tai-Peng Tsai$^1$}
\date{July 5, 2006}
\maketitle

{\footnotesize \noindent
$^1$
Department of Mathematics, University of British Columbia,
Vancouver, BC, Canada V6T 1Z2.
E-mail: gustaf@math.ubc.ca,  ttsai@math.ubc.ca\\
$^2$
Department of Mathematics,
Sungkyunkwan University and Institute of Basic Science,
Suwon 440-746, Republic of Korea.
E-mail: kkang@skku.edu}

\bigskip

\noindent{\bf Abstract:}\
We present new interior regularity criteria for suitable weak
solutions of the 3-D Navier-Stokes equations: a suitable weak solution
is regular near an interior point $z$ if either the scaled
$L^{p,q}_{x,t}$-norm of the velocity with $3/p+2/q\leq 2$, $1\leq
q\leq \infty$, or the $L^{p,q}_{x,t}$-norm of the vorticity with
$3/p+2/q\leq 3$, $1 \leq q < \infty$, or the
$L^{p,q}_{x,t}$-norm of the
gradient of the vorticity with $3/p+2/q\leq 4$, $1 \leq q$,
$1 \leq p$, is sufficiently small near $z$.


\section{Introduction}
We continue our study in \cite{GKT05} of the regularity problem for
{\it suitable weak solutions} $(u,p):\Om \times I \to \R^3 \times \R$
of the three-dimensional incompressible Navier-Stokes equations (NS)
\begin{equation}
\label{intro:nse}
\left\{
\begin{array}{c}
u_t-\Delta u+(u\cdot\nb)u+\nb p=f\\
{\rm div}\,\, u=0
\end{array}
\right.
\quad \mbox{ in }\,\,\Om \times I.
\end{equation}
Here $\Om$ is either a domain in $\R^3$ or the 3-dimensional torus
$\mathbb{T}^3$, $I$ is a finite time interval, $u(x,t)$ is the
velocity field and $p(x,t)$ is the pressure.  We also denote the
vorticity field by $w = \text{curl}\ u$. By suitable weak solutions we
mean functions which solve \eqref{intro:nse} in the sense of
distributions and satisfy some integrability conditions and the local
energy inequality (for details, see Definition \ref{defsws:100} in
section 2). For a point $z=(x,t)\in \R^3 \times \R$ we denote
\begin{equation*}
B_{x,r}= \{ y\in \R^3: |y-x|<r \},\quad Q_{z,r} :=B_{x,r} \times
(t-r^2, t).
\end{equation*}
A solution $u$ is said to be {\it regular at} $z \in \Om \times I$ if
$u\in L^\infty(Q_{z,r})$ for some $Q_{z,r}\subset \Om \times I$,
$r>0$. Otherwise it is {\it singular} at $z$ (see \cite[p.~780]{CKN}).

Although the existence of weak solutions was proved by Leray and Hopf
\cite{Leray,Hopf} in $\R^3$ and domains, it is not known whether the
solution stays regular for all time even if all the data are smooth.
One type of condition ensuring regularity involves
zero-dimensional integrals,
\begin{equation} \label{eq:1-5}
 \norm{u}_{L^{p,q}(\Om \times I)} < \infty, \quad \frac 3p + \frac 2q = 1,
 \quad 3\le p \le \infty,
\end{equation}
where
\begin{equation}
   \norm{u}_{L^{p,q}(\Om \times I)} =
\norm{u}_{L^q_t L^p_x(\Om \times I)} =
   \big \| \norm{u(x,t)}_{L^p_x(\Om)} \big \|_{L^q_t(I)}.
\end{equation}
These integrals have zero dimension if one assigns
the dimensions $1$, $2$, and $-1$ to $x$, $t$ and $u$. This is related
to the scaling property of solutions of (NS): The map
\begin{equation}
\label{eq:scaling}
\{u(x,t),p(x,t) \} \to \{\la u( \la x,\la ^2t),
\la^2p(\la x,\la ^2t)\}
\quad (\la>0),
\end{equation}
sends a solution of (NS) to another solution, with a new force
$\lambda^3 f(\lambda x, \lambda^2 t)$.

The first contributions in this direction,
concerning uniqueness and regularity of weak solutions,
were made by \cite{P59,S62,S63,L67} when
$3/p+2/q < 1$. The borderline cases $3/p+2/q = 1$, $3<p\le
\infty$, for different types of domains were later proved by
\cite{FJR, Sohr83, G86,S88}. See \cite{T90,Sohr01,CP} for results
in the setting of Lorentz spaces. The endpoint case
$(p,q)=(3,\infty)$ was recently resolved \cite{ESS} (also see the
references in \cite{Sohr01,ESS} for earlier results in
subclasses). Similar regularity criteria have been established
near the boundary \cite{VS02,K04,S05}.

In a series of papers \cite{VS76}--\cite{VS82}, Scheffer began to
study the partial regularity theory for (NS). His results were
further generalized and strengthened in Caffarelli-Kohn-Nirenberg
\cite{CKN}, which proved that the set $\calS$ of possible interior
singular points of a suitable weak solution is of one-dimensional
parabolic Hausdorff measure zero, i.e.  $\calP^1(\calS)=0$ (the
estimate of the Hausdorff measure was improved by a logarithmic
factor in \cite{CL00}).  The key to the analysis in \cite{CKN} is
the following regularity criterion: there is an absolute constant
$\eps >0$ such that, if $u$ is a suitable weak solution of (NS) in
$\Om \times I$ and if for an interior point $z\in \Om \times I$,
\begin{equation}
\limsup_{r\to 0_+}\frac{1}{r}\int_{Q_{z,r}}\abs{\nb u(y,s)}^2dyds
\leq \eps,
\label{nse:abs10}
\end{equation}
then $u$ is regular at $z$. See \cite{L98} for a simpler proof and
\cite{LS99} for more details.  See \cite{S02,SSS06} for extensions
when $z$ lies on a flat or curved boundary.

The objective of this paper is to present new sufficient conditions
for the regularity of suitable weak solutions to (NS) in the interior,
in terms of the smallness of the scaled $L^{p,q}$-norm of the
velocity, vorticity or the gradient of the vorticity.  We obtained
such results in terms of the velocity either in the interior or on a
flat boundary in \cite{GKT05}. We will assume that
the force $f$ belongs to a parabolic Morrey space $M_{2,\ga}$, for
some $ \ga > 0 $, equipped with the norm
\begin{equation}
  \norm{f}_{M_{2,\ga}(\Om \times I)}^2
= \sup_{Q_{z,r} \subset \Om \times I , \, r>0} \,
  \frac 1{r^{1+ 2 \ga}}  \int_{Q_{z,r}}\abs{f}^2dz'.
  \label{morrey:100}
\end{equation}
(This space is trivial if $\ga > 2$.)

Suitable weak solutions will be defined in Definition
\ref{defsws:100} of section 2.


\begin{thm}[Regularity Criteria] \label{TH1}
Suppose the pair $(u,p)$ is a suitable weak solution of (NS) in $\Om
\times I$ with force $f\in M_{2,\ga}(\Om \times I)$ for some
$\ga > 0$. Suppose $z=(x,t) \in \Om \times I$ and $Q_{z,r}
\subset \Om \times I$. Then $u$ is regular at $z$ if one of the
following conditions holds, for a small constant $\eps>0$ depending
only on $p^*$ (or $p,p^\sharp$), $q$, and $\ga$ (but independent of
$\norm{f}_{M_{2,\ga}}$).

(i) {\bf (Velocity criteria)} $u \in L^{p^*,q}_{\mathrm{loc}} $ near
$z$ and
\begin{equation}
\limsup_{r\to 0_+}\,\, r^{-( \frac 3{p^*} + \frac 2q -1)}
\norm{u-(u)_r}_{L^{p^*,q}(Q_{z,r})} \leq \eps,
\label{nse:abs30}
\end{equation}
where $(u)_r(s) = \frac 1{|B_r|} \int_{B_r} u(y,s)dy$,
for some $p^*,q$ satisfying
\begin{equation} \label{TH1:pq1}
  1\le 3/p^*+2/q \le 2, \quad 1\le p^*, q \le \infty.
\end{equation}

The same result holds if $u-(u)_r$ is replaced by $u$ in \eqref{nse:abs30}.

(ii) {\bf (Velocity gradient criteria)} $\nabla u \in
L^{p,q}_{\mathrm{loc}} $ near $z$ and
\begin{equation}
\limsup_{r\to 0_+}\,\, r^{-( \frac 3p + \frac 2q -2)}
\norm{\nabla u}_{L^{p,q}(Q_{z,r})} \leq \eps,
\label{nse:abs40}
\end{equation}
for some $p,q$ satisfying
\begin{equation} \label{TH1:pq2}
  2\le 3/p+2/q \le 3, \quad 1 \le q \le \infty.
\end{equation}

(iii) {\bf (Vorticity criteria)} $w= \curl u \in
L^{p,q}_{\mathrm{loc}} $ near $z$ and
\begin{equation}
\label{TH1:eq5}
\limsup_{r\to 0_+}\,\, r^{-( \frac 3p + \frac 2q -2)}
\norm{w}_{L^{p,q}(Q_{z,r})} \leq \eps,
\end{equation}
for some $p,q$ satisfying
\begin{equation} \label{TH1:pq3}
2\le 3/p+2/q \le 3, \quad 1 \le q \le \infty, \quad (p,q)\not = (1,\infty).
\end{equation}

(iv) {\bf (Vorticity gradient criteria)} $\nabla^2 u \in
L^{p^\sharp ,q}_{\mathrm{loc}} $ near $z$ and
\begin{equation}
\label{grad-vorticity}
\limsup_{r\to 0_+}\,\, r^{-(\frac 3{p^\sharp} + \frac 2q-3)}
\norm{\nabla w}_{L^{p^\sharp,q}(Q_{z,r})} \leq \eps,
\end{equation}
for some $p^\sharp, q$ satisfying
\begin{equation}\label{TH1:pq4}
3\le 3/p^\sharp + 2/q \le 4, \quad 1 \leq q, \quad 1 \le p^\sharp.
\end{equation}
Furthermore, for $p^\sharp > 1$,
$\nb w$ can be replaced by $\curl w$.

\end{thm}

{ \setlength{\unitlength}{2mm} \noindent
\begin{center}
\begin{picture}(30,28)
\put(1,0){\vector(1,0){26}}  \put (1,0){\vector(0,1){26}}
\put(1,24){\line(1,0){24}}   \put (25,0){\line(0,1){24}}
\put(9,0){\line(-2,3){8}}    \put (25,0){\line(-2,3){16}}
\put(17,0){\line(-2,3){16}}  \put (25,12){\line(-2,3){8}}
\put(1,24){\line(1,0){24}}  
\put(1,12){\line(1,0){8}}
\put(28,0){\makebox(0,0)[c]{\scriptsize $\frac 1p$}}
\put(0,26){\makebox(0,0)[c]{\scriptsize $\frac 1q$}}
\put(1,-1){\makebox(0,0)[c]{\scriptsize $0$}}
\put(9,-1){\makebox(0,0)[c]{\scriptsize $1/3$}}
\put(11.5,12){\makebox(0,0)[c]{\scriptsize $(\frac 13, \frac 12)$}}
\put(19.5,25){\makebox(0,0)[c]{\scriptsize $(\frac 23, 1)$}}
\put(27.5,12){\makebox(0,0)[c]{\scriptsize $(1, \frac 12)$}}
\put(1,24){\makebox(0,0)[c]{\scriptsize $\bullet$}}
\put(17,0){\makebox(0,0)[c]{\scriptsize$\bullet$}}
\put(17,24){\makebox(0,0)[c]{\scriptsize $\bullet$}}
\put(25,0){\makebox(0,0)[c]{\scriptsize o}}
\put(9,24){\makebox(0,0)[c]{\scriptsize $\bullet$}}
\put(9,12){\makebox(0,0)[c]{\scriptsize $\bullet$}}
\put(25,12){\makebox(0,0)[c]{\scriptsize $\bullet$}}
\put(11,25){\makebox(0,0)[c]{\scriptsize $(\frac 13, 1)$}}
\put(17,-1){\makebox(0,0)[c]{\scriptsize $2/3$}}
\put(25,-1){\makebox(0,0)[c]{\scriptsize $1$}}
\put(0,0){\makebox(0,0)[c]{\scriptsize $0$}}
\put(0,12){\makebox(0,0)[c]{\scriptsize $\frac 12$}}
\put(0,24){\makebox(0,0)[c]{\scriptsize $1$}}
\put(3,4){\makebox(0,0)[c]{\scriptsize I}}
\put(10,4){\makebox(0,0)[c]{\scriptsize II}}
\put(4,16){\makebox(0,0)[c]{\scriptsize III}}
\put(11,16){\makebox(0,0)[c]{\scriptsize IV}}
\put(18,16){\makebox(0,0)[c]{\scriptsize V}}
\put(13,-3){\makebox(0,0)[c]{\scriptsize Figure 1: Regularity
Criteria}}
\end{picture}
\end{center}
}
%

\bigskip

{\it Comments for Theorem \ref{TH1}.}

\begin{enumerate}

\item The region defined by \eqref{TH1:pq1} corresponds to the union
of II and III in Figure 1, including all borderlines. The region
defined by \eqref{TH1:pq2} corresponds to IV, including all
borderlines.  The region defined by \eqref{TH1:pq3} also corresponds
to IV, but without the corner point
$(1/p,1/q)=(1,0)$.  The region defined by \eqref{TH1:pq4} corresponds
to V.

\item In \eqref{TH1:pq1}, the lower bound $1 \le 3/p^* + 2/q$ is only
to ensure a non-positive exponent of $r$ in \eqref{nse:abs30}. The
true limit is the upper bound $3/p^* + 2/q \le 2$. Similar comments
apply to \eqref{TH1:pq2}, \eqref{TH1:pq3} and \eqref{TH1:pq4}.

\item The quantities in \eqref{nse:abs30}, \eqref{nse:abs40},
\eqref{TH1:eq5} and \eqref{grad-vorticity} are zero-dimensional, and are
invariant under the scaling \eqref{eq:scaling}.  Such quantities are
useful in the regularity theory for (NS), see e.g.~\cite{CKN}.

\item In \cite{GKT05}, the authors obtained Theorem \ref{TH1} (i) only
for region II, without the borderline $q=2$ (but the result is
also valid on a flat boundary of $\Om$). Theorem \ref{TH1} (i)
extends it to region III, and in particular includes the point
$(1/p,1/q)=(0,1/2)$.  It does not further assume the smallness of
the pressure, in contrast to, e.g., Theorem \ref{mod:lemma}.
Special cases $(1/p,1/q)=(1/3,1/3)$ and $(1/2,0)$ were obtained in
\cite{TX} and \cite{SS}, respectively.

\item Theorem \ref{TH1} (ii) contains the special case $(p,q)=(2,2)$
of \cite{CKN}.

\item Theorem \ref{TH1} (iii) contains the special case $(p,q)=(2,2)$
of \cite{TX}.

\end{enumerate}

Theorem \ref{TH1} implies many known regularity criteria. Some of
them are summarized below. For simplicity we assume $f=0$. The
Lorentz space $L^{(p,\infty)}$ for $p<\infty$ is defined with the
norm $\norm{v}_{L^{(p,\infty)}} = \sup _{\si>0} \si
|\{|v|>\si\}|^{1/p}$.

\begin{cor}\label{TH4}
Let $u$ be a weak solution of (NS) in $\Om \times I$
with $f=0$ and
$Q_{z_0,r_0}\subset \Om \times I$ for some $r_0>0$. Then $u$ is regular
at $z_0$ if one of the following conditions holds.

(i) (zero-dimensional integrals of $u$ \cite{FJR, Sohr83, G86,S88})
If
\begin{equation}
u\in L^{p,q}(Q_{z_0,r_0}), \qquad
\frac{3}{p}+\frac{2}{q}= 1,\quad 3<p\leq \infty, \label{SPC:sec3}
\end{equation}
or $u\in L^{3,\infty}(Q_{z_0,r_0})$ and
$\norm{u}_{L^{3,\infty}(Q_{z_0,r_0})}$ is sufficiently
small.

(ii) (Lorentz spaces \cite{T90,Ko98,Sohr01,CP}) If $u$ is in
$L^{(q,\infty)}((t_0-r^2,t_0); L^{(p,\infty)}(B_{x_0,r}))$ with
$3/p+2/q=1$, $3 < p < \infty$, and
$\norm{u}_{L^{(q,\infty)}_t\,L^{(p,\infty)}_x(Q_{z_0,r})}$ is
sufficiently small.

(iii) (zero-dimensional integrals of $\nabla u$ \cite{BdV95})
\begin{equation*}
\nb u\in L^{p,q}(Q_{z_0,r_0}), \qquad
\frac{3}{p}+\frac{2}{q}=2,\quad\frac{3}{2} < p \le \infty,
\end{equation*}
or $\nb u\in L^{3/2,\infty}(Q_{z_0,r_0})$ and
$\norm{\nb u}_{L^{3/2,\infty}(Q_{z_0,r_0})}$
is sufficiently small.

(iv) (zero-dimensional integrals of $w= \curl u$ \cite{CKL})
\begin{equation}\label{vor:sec3}
w\in L^{p,q}(Q_{z_0,r_0}), \qquad
\frac{3}{p}+\frac{2}{q}=2,\quad\frac{3}{2} < p \leq \infty,
\end{equation}
or $w \in L^{3/2,\infty}(Q_{z_0,r_0})$ and
$\norm{w}_{L^{3/2,\infty}(Q_{z_0,r_0})}$
is sufficiently small.

\end{cor}

{\it Comments for Corollary \ref{TH4}.}

\begin{enumerate}

\item To prove Corollary \ref{TH4} using Theorem \ref{TH1},
we need to show that $u$ is suitable under the corresponding
assumptions.  It suffices to show that $|u|^2 |\nabla u| \in
L^1_{t,x}$, which justifies the integration by parts and thus one can
prove the local energy inequality. In fact, it is enough to show $u
\in L^4_{t,x}$ since $\iint |u|^2 |\nabla u| dz \leq \norm{u}_{L^4}^2
\norm{\nb u}_{L^2}$.

For (i), it follows from $\norm{u}_{L^4}^2\le \norm{u}_{L^{p,q}}
\norm{u}^{2/q}_{L^{2,\infty}} \norm{u}^{3/p}_{L^{6,2}}$.

For (ii), since $3<p<\infty$, one can choose $p_1$,$q_1$ so that
\[
q_1 < q,\quad p _1 < p, \quad 1/p_1 + 1/q_1 \le 1/2, \quad
 3/p_1 + 1/q_1 \le 1.
\]
That is, $(1/p_1,1/q_1)$ lies in region V of Figure 2 of \cite{GKT05}.
By the imbedding of $L^{(p,\infty)} \subset L^{p_1}$ and $L^{(q,\infty)}
\subset L^{q_1}$, we have $u \in L^{p_1,q_1}$.  Interpolating with $u
\in L^{2,\infty} \cap L^{6,2}$, we get $u \in L^4_{t,x}$.

For (iii), we have $\displaystyle \iint |u|^2 |\nabla u| dz \leq
\|u\|^{2/q}_{L^{2,\infty}_{x,t}} \| u\|^{3/p}_{L^{6,2}_{x,t}}
\|\nabla u \|_{L^{p,q}_{x,t}}$.

For (iv), since $\|\nabla u\|_{L^{p,q}(Q_r)}\leq
C\|w\|_{L^{p,q}(Q_{2r})}+ C\|u\|_{L^{p,q}(Q_{2r})}$ (see Remark
\ref{RK3.7}), it follows from (iv).

\item Strictly speaking, one also needs to show that $p \in L^{3/2}$
so that $(u,p)$ is suitable. But this has already been done
\cite{SW86,L98}. By \cite[Lem.~3.4]{L98}, one has $\nb p \in L^{5/3}_t
L^{15/14}_x(Q_r)$ for every week solution in $Q_r$. Let $\tilde p(x,t)
= p(x,t) - \aint_{B_r} p(x,t)\, dx$.  The new pair $(u,\tilde p)$ is
suitable since the local energy inequality \eqref{lei} remains the
same if one replaces $p$ by $\tilde p$, and $\tilde p \in
L^{5/3}_{t,x}(Q_r)$ by Poincar\'e inequality.

\item We now complete the proof of Corollary \ref{TH4}.  For (ii),
since $3<p<\infty$, one can choose $q_2 < q$, $p _2 < p$, and $3/p_2 +
2/q_2 = 2$.  Being small in $L^{(q,\infty)}_t L^{(p,\infty)}_x(Q_{z_0,r})$ implies
smallness in the scaled norm $\frac 1r L^{q_2}L^{p_2}(Q_r)$ by
imbedding. Then one applies Theorem \ref{TH1}.  For the rest, one
imbeds $L^{p,q}$ to $L^{p_2,q}$ for some suitable $p_2 < p$.

\item Corollary \ref{TH4} (i) is due to several authors, already
quoted above. Theorem \ref{TH1} does not imply the end point case
$u\in L^{3,\infty}(Q_{z_0,r_0})$ without smallness assumption, for
which see \cite{ESS}.

\item For Corollary \ref{TH4} (ii), \cite{T90} proved regularity for
small $u$ in the classes $L^{(q,\infty)}_tL^p_x(Q_r) $ with $3<p<\infty$.
\cite{Ko98} in the class $L^{\infty}_tL^{(3,\infty)}_x$
(see \cite{Ko01} for improvement).
\cite{Sohr01} in the classes $L^{(q,\infty)}_tL^{(p,\infty)}_x(\Om \times
I) $ with $3<p<\infty$.
\cite{CP} in the classes $L^{(q,\infty)}_tL^{(p,\infty)}_x $ with
$3<p<\infty$ and the classes $L^{(p,\infty)}_x L^{(q,\infty)}_t$ with
$3\le p<\infty$.
It follows from these results, in particular, that $u$ is regular at
$z_0$ if it satisfies, for $\theta \in [0,1]$ and some
$\eps=\eps(\theta) >0$,
\begin{equation}\label{TH4:eq2}
\lim_{r\to 0}{\rm ess}\, \sup_{Q_{z_0,r}} \abs{t-t_0}^{\theta/2}
\abs{x-x_0}^{1-\theta} \abs{u(x,t)}\leq \eps.
\end{equation}
Our Theorem \ref{TH1} does not cover the endpoint cases $p=3,\infty$,
except the cases $\th=0,1$ in \eqref{TH4:eq2} when suitability is
assumed.

\item Corollary \ref{TH4} (iii) was proved in \cite{BdV95} for the
cases $3/2<p<\infty$.
The endpoints $p=3/2$ and $p=\infty$ were not
obtained in \cite{BdV95}.
The $p=3/2$ case without the smallness assumption follows from
\cite{ESS} and imbedding.

\item Corollary \ref{TH4} (iv) was proved in \cite[Prop.~2]{CKL}.  The
main result in \cite[Th.~1]{CKL} shows regularity near $z_0$ assuming
only two components of the vorticity belonging to $L^{p,q}_{x,t}$.
Again, the $p=3/2$ case without the smallness
assumption follows from~\cite{ESS} and Remark \ref{RK3.7}.

\end{enumerate}

A major motivation for the study of such regularity criteria is to
improve the partial regularity result of \cite{CKN}.  For example,
Constantin \cite{Co90} proved, when $\Om=\mathbb{T}^3$, the existence
of suitable weak solutions satisfying
\begin{equation}
\label{Const}
  \nabla w \in L^{4/3-\e}(\Om \times I), \quad \forall 0<\e \ll 1.
\end{equation}
Note that the integral $\iint |\nabla w|^{4/3-\e} dz$ has dimension
$1+3\e$.  Combining this estimate with Theorem \ref{TH1} (iv), we find
that the parabolic Hausdorff dimension of the singular set
$\mathcal{S}$ of $u$ is at most one. This is slightly weaker than the
CKN theorem that the one-dimensional parabolic Hausdorff measure of
$\mathcal{S}$ is actually zero. Note that Scheffer \cite{VS85,VS87}
constructed examples satisfying the local energy inequality and their
dimensions of singular sets are arbitrarily close to one. Thus the CKN
result is optimal for functions satisfying only the local energy
inequality.  However, the proof of \eqref{Const} uses the equation for
the vorticity, which may not be satisfied by Scheffer's
examples. Therefore there might be hope to prove other a priori
estimates for $w$ and thus improve the partial regularity.

The rest of this paper is organized as follows. In Section 2 we
introduce some scaling invariant functionals, recall the notion of
suitable weak solutions and a regularity criterion involving the
scaled norms of velocity and pressure. In Section 3 we establish some
estimates regarding the velocity, pressure and vorticity, and prove
Theorem \ref{TH1}.


\section{Preliminaries}

In this section we introduce the notation, review suitable weak
solutions, and recall a regularity criterion involving scaled norms.

We start with the notation.  Let $\Om$ be either an open domain in
$\R^3$ or the 3-dimensional torus $\mathbb{T}^3$, and $I$ be a finite
time interval.  By $N=N(\alpha,\beta,\ldots)$ we denote a constant depending
on the prescribed quantities $\alpha,\beta,\ldots$, which may change
from line to line. For $1\leq q\leq\infty$, $W^{k,q}(\Om)$ denote the
usual Sobolev spaces, i.e. $W^{k,q}(\Om)= \{ f\in
L^q(\Om):D^{\alpha}f\in L^q(\Om), 0\leq|\alpha|\leq k\}$.
We denote by $\aint_E f$ the average of $f$ on $E$; i.e., $\aint_E
f=\int_E f/|E|$. For a
function $f(x,t)$, $E \subset \Om$ and $J \subset I$, we denote
$
\norm{f}_{L^{p,q}(E \times J)} =
\norm{f}_{L^q L^p(E \times J)} = \big\| \norm{f}_{L^p(E)} \big \|_{L^q(J)} $.

Next, we define several scaling-invariant functionals similar to
those in \cite{CKN, L98, LS99, S02}. For a suitable weak solution
$(u,p)$ and $z=(x,t)\in\Om\times I$, let
\begin{equation*}
A(r):=\sup_{t-r^2\leq s<t}
\frac{1}{r}\int_{B_{x,r}}\abs{u(y,s)}^2dy,\quad
E(r):=\frac{1}{r}\int_{Q_{z,r}}\abs{\nb u(y,s)}^2\ind y \ind s,
\end{equation*}
\begin{equation*}
C(r):=\frac{1}{r^2}\int_{Q_{z,r}}\abs{u(y,s)}^3\ind y\ind s,
\quad
\tilde{C}(r):=\frac{1}{r^2}\int_{Q_{z,r}}\abs{u(y,s)-(u)_r(s)}^3\ind
y\ind s,\quad
\end{equation*}
\begin{equation*}
D(r):= \frac{1}{r^2} \int_{Q_{z,r}}
\abs{p(y,s)}^{\frac{3}{2}} \ind y \ind s.
\end{equation*}
where $(u)_r(s)=\frac 1{|B_{x,r}|}\int_{B_{x,r}}u(\cdot,s)dy$.
Let $p,q$ and $p^*$ be numbers satisfying
\begin{equation} \label{pq-100}
\frac{3}{p}+\frac{2}{q}=3,\quad 1\leq q\leq \infty,\quad
\frac{1}{p^*}=\frac{1}{p}-\frac{1}{3}.
\end{equation}
Recall $w=\nabla\times u$ is the vorticity field of $u$. We define
\begin{equation*}
\tilde G(r):= \frac 1r \norm{u(y,s)-(u)_r(s)}_{L^q_sL^{p^*}_y(Q_{z,r})},
\quad
G_1(r):= \frac 1r \norm{\nabla u(y,s)}_{L^q_sL^p_y(Q_{z,r})},
\end{equation*}
\[
W(r):= \frac 1r \norm{w(y,s)}_{L^q_sL^p_y(Q_{z,r})}.
\]
When $1\le q \le 2$, we also define
\[
W_1(r):= \frac 1r \norm{\nabla w(y,s)}_{L^q_sL^{p^\sharp}_y(Q_{z,r})}, \quad
\tilde{W}_1(r) := \frac{1}{r}\norm{\curl w(y,s)}_{L^q_sL^{p^\sharp}_y(Q_{z,r})},\end{equation*}
where  $p^\sharp$
is the number satisfying, for $p,q$ as in \eqref{pq-100},
\begin{equation}\label{grad-vorticity100}
\frac 3{p^\sharp} + \frac 2q = 4,\quad \frac 1{p^{}} = \frac
1{p^\sharp} - \frac 13,\quad 1\leq p^\sharp\leq \frac{3}{2}.
\end{equation}

\medskip

We now define suitable weak solutions for the (NS).

\begin{defi}\label{defsws:100}
Suppose that $f$ belongs
to the parabolic Morrey space $M_{2,\ga}(\Om \times I)$ for some
$\ga \in (0,2]$.  A pair $(u,p)$ is a {\bf suitable weak solution}
to the Navier-Stokes equations (\ref{intro:nse}) in $\Om \times I$ with
force $f$ if the following conditions are satisfied.
\begin{itemize}
\item[(a)]
The functions $u:\Om \times I \to \R^3$ and $p:\Om \times I \to \R$ satisfy
\begin{equation}
u\in L^{\infty}(I;L^2(\Om))\cap L^2(I;W^{1,2}(\Om)),\quad p\in
L^{\frac{3}{2}}(\Om \times I). \label{defsws:200}
\end{equation}

\item[(b)] $u$ and $p$ solve (\ref{intro:nse}) in $\Om \times I$ in the
sense of distributions.

\item[(c)]
$u$ and $p$ satisfy the local energy inequality

\parbox{11cm}{
\[
\int_{\Om}\abs{u(x,t)}^2\phi(x,t)\ind x +2\int_{t_0}^t \int_\Om \abs{\nb
u(x,t')}^2\phi(x,t')\ind x\ind t'
\]
\[
\leq \int_{t_0}^t \int_\Om \bke{\abs{u}^2(\pr_t\phi+\Delta\phi)
+(\abs{u}^2+2p)u\cdot\nb\phi
+2f\cdot u\phi}\ind x\ind t'
\]} \hfill
\parbox{1cm}{\begin{equation}\label{lei}\end{equation}}
for all $t\in I =(t_0,t_1) $ and all nonnegative functions
$\phi\in\calC^{\infty}_0(\Om \times I)$. \qed
\end{itemize}
\end{defi}

In this definition we impose no initial or boundary condition for
$u$.

The main difference between suitable weak solutions and {\it
Leray-Hopf weak solutions} (see \cite[p.779]{CKN}) is the additional
condition of the local energy inequality \eqref{lei}.
The existence of suitable weak solutions is proved in \cite{VS77,CKN}.
Definition \ref{defsws:100} is the slightly modified version used in
\cite{L98}. As remarked in \cite[page 823]{CKN}, it is an open
question if all weak solutions are suitable.

Next we recall a local regularity criterion, which is a refined
version of \cite[Prop.~1]{CKN}, and is formulated in the present form
with $f=0$ in \cite{NRS96,L98}, and proved with nonzero $f\in
M_{2,\ga}$ in \cite[Prop.~2.8]{LS99}.

\begin{thm}\label{mod:lemma}
There exists $\eps > 0$ depending only on $\ga>0$ (and independent
of $\norm{f}_{M_{2,\ga}}$), such that if $(u,p)$ is a suitable weak
solution of (NS) with $f \in M_{2,\ga}$, then $u$ is regular at
$z=(x,t)\in\Om\times I$ if
\begin{equation}\label{CKN-cond}
 C(r)+D(r)<\eps \quad \text{for some }r>0.
\end{equation}
\end{thm}

An important feature of \eqref{CKN-cond} is that it requires only one
$r$, not infinitely many $r$. We will prove our regularity criteria
based on this theorem.  For our proof in the next section, in order to get
\eqref{CKN-cond}, it suffices to assume $\ga > -1$. The assumption
$\ga >0$ is made in order to apply Theorem \ref{mod:lemma}.

\section{Local interior regularity}
In this section, we present the proof of Theorem \ref{TH1}. Through
the entire section, we assume $(u,p)$ is a suitable weak solution in
$\Om \times I$. Without loss of generality, we assume $z = (0,0)$ and
$Q_r= Q_{(0,0),r} \subset \Om \times I$. By H\"older inequality, it
suffices to consider borderline exponents, i.e., those exponents $p,
p^*,p^\sharp$ and $q$ satisfying \eqref{pq-100} and
\eqref{grad-vorticity100}.
Denote $m_\ga = \norm{f}_{M_{2,\ga}}$.

\begin{lemm}\label{lei2}
Suppose $Q_{2r} \subset \Om \times I$ and $0< r \le
m_\ga
^{-1/(1+\ga)}$.  Then
\[
A(r)+E(r) \leq N[1+C(2r)+D(2r)].
\]
\end{lemm}
\begin{pf}
By choosing suitably localized $\phi$ in the local energy
inequality \eqref{lei}, we get
\[
A(r)+E(r) \leq N\left(C^{\frac{2}{3}}(2r)+C(2r)
+\frac{1}{r^2}\norm{u}_{L^{3}(Q_{2r})}\norm{p}_{L^{\frac{3}{2}}(Q_{2r})}
+r\int_{Q_{2r}}|f|^2dz'\right)
\]
which is bounded by $ N[1+C(2r)+D(2r) +r^{2(\ga+1)}m^2_{\ga}]$.
\end{pf}

\begin{lemm}\label{basiclemma}
Suppose $u\in L^{p^*,q}(Q_r)$ with $3/p^*+2/q=2$, $1 \le q \le
\infty$, then
\begin{equation*}
\tilde C(r)\leq N A^{\frac{1}{q}}(r) E^{1-\frac{1}{q}}(r)\tilde G(r).
\end{equation*}
\end{lemm}
\begin{pf}
Let $\alpha=(2p^*-3)/3p^*$ and $\beta=1/p^*$. Note
$1/3=\alpha/2+\beta/6+(1-\alpha-\beta)/p^*$. Using the H\"older
inequality and Sobolev imbedding, we obtain
\[
\begin{split}
\norm{u-(u)_{r}}_{L^3(B_r)} &\leq
N\norm{u}^{\alpha}_{L^2(B_r)}\norm{u-(u)_{r}}^{\beta}_{L^6(B_r)}
\norm{u-(u)_{r}}^{1-\alpha-\beta}_{L^{p^*}(B_r)} \\
&\leq N\norm{u}^{\alpha}_{L^2(B_r)}\norm{\nabla
u}^{\beta}_{L^2(B_r)}\norm{u-(u)_{r}}^{\frac{1}{3}}_{L^{p^*}(B_r)},
\end{split}
\]
where we used $1-\alpha-\beta=1/3$. Raising to the third power,
integrating in time and dividing both sides by $r^2$, we get
\[
\tilde C(r)\leq
\frac {N}{r^2}\int_{-r^2}^0\norm{u}^{3\alpha}_{L^2(B_r)}\norm{\nabla
u}^{3\beta}_{L^2(B_r)}\norm{u-(u)_{r}}_{L^{p^*}(B_r)}dt
\]
\[
\leq
\frac {N}{r^2}
r^{\frac{3}{2}\alpha}A^{\frac{3}{2}\alpha}(r)\bke{\int_{-r^2}^0\norm{\nabla
u}^2_{L^2(B_r)}dt}^{\frac{3\beta}{2}}\bke{\int_{-r^2}^0\norm{
u-(u)_{r}}^q_{L^{p^*}(B_r)}dt}^{\frac{1}{q}},
\]
which equals $N A^{\frac{1}{q}}(r)E^{1-\frac{1}{q}}(r)\tilde G(r)$.
\end{pf}

\begin{lemm}\label{L3-1}
Suppose $0 < 2r \leq \rho$ and $Q_{\rho}\subset \Om \times I$. Then
\begin{equation*}
C(r)\leq
N\left(\frac{r}{\rho}\right)C(\rho)+
N\left(\frac{\rho}{r}\right)^2\tilde{C}(\rho).
\end{equation*}
\end{lemm}
\begin{pf}
This follows from the H\"older inequality:
\[
C(r) \leq \frac N{r^2}  \int_{Q_r} \bke{|(u)_{\rho}|^3
+|u-(u)_{\rho}|^3} \,dz' \leq N(\frac r\rho) C(\rho) + N(\frac \rho
r)^2\tilde{C}(\rho).
\]
\end{pf}

\begin{lemm}\label{preest:100}
Suppose $0 < 2r \leq \rho$ and $Q_{\rho}\subset \Om \times I$. Then
\begin{equation}
D(r)\leq N\bke{\frac{\rho}{r}}^2
(\tilde C(\rho)+\rho^{\frac{3}{2}(\ga+1)}m^{\frac{3}{2}}_{\ga}) +
N\bke{\frac{r}{\rho}}D(\rho).
\label{pres:100}
\end{equation}
\end{lemm}
\begin{pf}  Let $\phi(x)\ge 0$ be supported in $B_{\rho}$ with
$\phi=1$ in $B_{\rho/2}$. The divergence of \eqref{intro:nse} gives
$-\Delta p=\partial_{i} \partial_{j}\bke{u_iu_j} - \nabla \cdot f$ in
the sense of distributions. Let
\[
p_1(x,t):=\int_{\R^3}\frac{4\pi}{\abs{x-y}}\bket{\partial_{i}
\partial_{j}\bkt{ (u_i-(u_i)_\rho)(u_j-(u_j)_\rho)\phi}
-\nabla\cdot(f\phi)}(y,t)dy
\]
and $p_2(x,t):=p(x,t)-p_1(x,t)$. Due to $\div u=0$, $\Delta
p_2 =0$ in $B_{\rho/2}$.  By the mean value property of harmonic
functions,
\[
\frac{1}{r^2}\int_{B_r}\abs{p_2}^{\frac{3}{2}}dx\leq
\frac{Nr}{{\rho}^3}\int_{B_{\rho/2}}\abs{p_2}^{\frac{3}{2}}dx
\leq
\frac{Nr}{{\rho}^3}\int_{B_{\rho}}|p|^{\frac32}dx
+ \frac{Nr}{{\rho}^3}\int_{B_{\rho}} \abs{p_1}^{\frac{3}{2}}dx.
\]
By Calderon-Zygmund and potential estimates,
\[
\frac r{\rho^3}\int_{B_\rho} \abs{p_1}^{\frac{3}{2}}dx\leq
\frac 1{r^2}\int_{B_\rho} \abs{p_1}^{\frac{3}{2}}dx\le
\frac
N{r^2}\int_{B_{\rho}}\abs{u-(u)_\rho}^3
+\frac{N\rho^{9/4}}{r^2}\big(\int_{B_{\rho}}|f|^2dx
\big)^{\frac{3}{4}}.
\]
Adding these estimates, integrating in time, and using $ \int
_{-r^2}^0 \frac {\rho^{9/4}}{r^{2}} \big(\int_{B_{\rho}}|f|^2dx
\big)^{\frac{3}{4}} dt\le N r^{-3/2}m_\ga^{3/2}\rho^{3+3\ga/2} $, we
get
\[
\frac{1}{r^2}\int_{Q_r}\abs{p}^{\frac{3}{2}}dz'\leq
\frac{1}{r^2}\int_{Q_r}\abs{p_1}^{\frac{3}{2}}+
\abs{p_2}^{\frac{3}{2}}dz'
\le \text{RHS of } \eqref{pres:100}.
\]
\end{pf}

Now we are ready to prove Theorem \ref{TH1} (i).

\MAINPF{(i)}

It suffices to prove the borderline cases
$3/p^*+2/q=2$ and $1 \le q \le \infty$.
The other cases follow by H\"older inequality.
Suppose $0<4r\le \rho$.  By Lemmas \ref{basiclemma}
and \ref{preest:100}, and by Lemma \ref{L3-1}, we get
\[
C(r)+D(r)\leq N(\frac{r}{\rho})\bke{C(\frac{\rho}{2})
+D(\frac{\rho}{2})}+N(\frac{\rho}{r})^2
\bke{\tilde{C}(\frac{\rho}{2})+\rho^{\frac32(\ga+1)}m^{\frac{3}{2}}_{\ga}}
\]
\[
\leq N(\frac{r}{\rho})\bke{C(\frac{\rho}{2})
+D(\frac{\rho}{2})}+N(\frac{\rho}{r})^2
\bke{A^{\frac{1}{q}}(\frac{\rho}{2})E^{1-\frac{1}{q}}(\frac{\rho}{2})\tilde G(\frac{\rho}{2})
+\rho^{\frac32(\ga+1)}m^{\frac{3}{2}}_{\ga}}.
\]
Suppose  $\rho\le
m_\ga^{-1/(\ga+1)}$. By Lemma \ref{lei2},
\[
N(\frac{\rho}{r})^2
A^{\frac{1}{q}}(\frac{\rho}{2})E^{1-\frac{1}{q}}(\frac{\rho}{2})\tilde G(\frac{\rho}{2})
\leq N(\frac{\rho}{r})^2 \bke{1+C(\rho)+D(\rho)}\tilde G(\rho).
\]
Combining the above estimates, we obtain
\[
C(r)+D(r)\leq
N_2\bke{(\frac{r}{\rho})+(\frac{\rho}{r})^2\tilde G(\rho)}\bke{C(\rho)
+D(\rho)}
+N_2(\frac{\rho}{r})^2\bke{\tilde G(\rho)
+\rho^{\frac32(\ga+1)}m^{\frac{3}{2}}_{\ga}}.
\]
Choose $\theta\in (0,1/4)$ so that $N_2\theta<1/4$.  We fix $r_0<
\min \{1, \frac {1}{m_\ga},
\frac 1{m_\ga} (\frac{\eps \th^2}{8 N_2})^{2/3}\}^{1/(\ga+1)}$
such that $\tilde G(r)<\frac{\theta^2}{1+8N_2} \min \{1,\e\}$ for all $r\leq r_0$,
where $\eps$ is the constant in Theorem
\ref{mod:lemma}. Replacing $r$ and $\rho$ by $\theta r$ and $r$,
respectively, we get
\[
C(\theta r)+D(\theta r)\leq \frac{1}{2} \bke{C(r)+D(r)}+\frac \e4,
\quad \forall r<r_0.
\]
By iteration,
\begin{equation*}
C(\theta^k r)+D(\theta^k r) \leq
\frac{1}{2^k}\bke{C(r)+D(r)}+\frac{\eps}{2}, \quad \forall r<r_0.
\end{equation*}
Thus, for $k$ sufficiently large, $C(\theta^k r)+D(\theta^k r) \le
\e$, from which $z$ is a regular point due to Theorem \ref{mod:lemma}.

The last statement of Theorem \ref{TH1} (i), that one can replace
$u-(u)_r$ by $u$, is because $\norm{u - (u)_r}_{L^{p^*,q}} \le N
\norm{u}_{L^{p^*,q}} $.  \qed

\bigskip

The following modification of Lemma~\ref{basiclemma}
is all that is needed to prove Theorem~\ref{TH1} (ii).
\begin{lemm}\label{lemma3-6}
Suppose $0<2r\le \rho$ and $Q_{\rho}\subset \Om \times I$. Then
\begin{equation} \label{lemma3-6a}
\tilde{C}(r)  \leq N A^{1/q}(r)E^{1-1/q}(r)G_1(r),
\end{equation}
\end{lemm}
\begin{pf}
The proof is similar to that of Lemma \ref{basiclemma}. When $1\le
p<3$, using the same exponents $\alpha=1-1/p$ and $\beta=1/p -1/3$, we
have
\begin{align*}
\norm{u-(u)_{r}}_{L^3(B_r)}^3 &\leq N\norm{u}^{3
\alpha}_{L^2(B_r)}\norm{u-(u)_{r}}^{3\beta}_{L^6(B_r)}
\norm{u-(u)_{r}}^{3(1-\alpha-\beta)}_{L^{p^*}(B_r)}
\\
&\leq N\norm{u}^{2/q}_{L^2(B_r)}\norm{\nabla
u}^{2-2/q}_{L^2(B_r)}\norm{\nabla u}_{L^{p}(B_r)}.
\end{align*}
If $p=3$ (and $q=1$), by Gagliardo-Nirenberg and Poincar\'e inequalities,
\begin{align*}
\norm{u-(u)_{r}}_{L^3(B_r)}^3 &\leq N \norm{u-(u)_{r}}_{L^2(B_r)}^2
\norm{\nabla u}_{L^{3}(B_r)}+
 \frac N {r^{3/2}} \norm{u-(u)_{r}}_{L^2(B_r)}^3
\\
&\leq N \norm{u}_{L^2(B_r)}^2 \norm{\nabla u}_{L^{3}(B_r)}.
\end{align*}
Integrating in time and applying the H\"older inequality, we get
\eqref{lemma3-6a}.
\end{pf}

\MAINPF{(ii)}

The proof is the same as that for Theorem \ref{TH1} (i): we only need
to replace Lemma \ref{basiclemma} by Lemma \ref{lemma3-6}, and
replace the quantity $\tilde G(r)$ by $G_1(r)$.  \qed

\bigskip

The next lemma shows that the gradient of the velocity
can be controlled by the vorticity. This is the key
to Theorem~\ref{TH1} (iii).
\begin{lemm}\label{TH3-4}
Suppose $0<2r\le\rho$ and $Q_{\rho} \subset \Om \times I$.
Suppose $\nabla u\in L^{p,q}_{x,t}(Q_{\rho})$ with
$\frac 3p + \frac 2q =3$ and $1 \le q < \infty$. Then
\begin{equation}
  G_1(r)\leq N\bke{\frac{\rho}{r}}W(\rho)
  +N\bke{\frac{r}{\rho}}^{\frac{3}{p}-1}G_1(\rho).
\label{wv-100}
\end{equation}
Furthermore, if $p=3$ (so $q=1$), then
\begin{equation}
  G_1(r) \leq N\bke{ \frac{\rho}{r}}W(\rho)
  + N\bke{\frac{r}{\rho}}G_1(\rho)
  + g(u;r)
\label{wv-101}
\end{equation}
where $g(u;r) \to 0$ as $r \to 0$.
\end{lemm}

\begin{pf}
Choose a standard cut off function $\phi$ supported in
$B_{\rho}$ such that $\phi=1$ in $B_{3\rho/4}$. Define
\[
v(x,t):=\int_{\R^3}\nabla_x \frac{4\pi}{|x-y|} \times w(y,t)\phi(y)dy,
\quad h=u-v.
\]
Note that $\Delta_x h(x,t) = 0$ in $B_{3\rho/4}$.

We give the proof of~\eqref{wv-100} first.
By the mean value
property of harmonic functions, for each fixed time $t$,
\begin{equation*}
\norm{\nb h}_{L^p(B_r)} \le N(\frac{r}{\rho})^{3/p} \norm{\nb
h}_{L^p(B_{\rho/2})} \le N(\frac{r}{\rho})^{3/p} \bke{\norm{\nb
u}_{L^p(B_{\rho})} + \norm{\nb v}_{L^p}}.
\end{equation*}
On the other hand, due to Calderon-Zygmund estimates,
for each fixed time,
\[
\norm{\nb v}_{L^p} \le N \norm{w}_{L^p(B_\rho)}.
\]
Combining these estimates, we obtain
\[
\norm{\nabla u}_{L^{p}(B_r)}\leq \norm{\nabla
v}_{L^{p}(B_r)}+\norm{\nabla h}_{L^{p}(B_r)}
\leq N \norm{w}_{L^{p}(B_{\rho})}
+N (\frac{r}{\rho})^{\frac{3}{p}}\norm{\nabla
u}_{L^{p}(B_{\rho})}.
\]
Taking $L^q$-norm in time and dividing both sides by $r$, we get
\eqref{wv-100}.

To prove~\eqref{wv-101}, set $p=3$ (so $q=1$),
use the above estimate for $\nb v$, and
modify the estimate for $\nb h$ as follows:
\begin{equation}
\label{eq:harm}
  \norm{\nb h}_{L^3(B_r)}
  \le \norm{\nb h - (\nb h)_r}_{L^3(B_r)}
  + \norm{(\nb h)_r}_{L^3(B_r)}.
\end{equation}
The second term in~\eqref{eq:harm} is just
$N r |(\nb h)_r|$.
For the first term in~\eqref{eq:harm},
use the Poincar\'e-Sobolev inequality,
the mean-value property, and an interior estimate:
\begin{align}  
\nonumber
  \norm{\nb h - (\nb h)_r}_{L^3(B_r)} &
  \le N \norm{\nb^2 h}_{L^{3/2}(B_r)}
  \le N (\frac{r}{\rho})^2 \norm{\nb^2 h}_{L^{3/2}(B_{\rho/2})} \\
\nonumber
  &\le N (\frac{r}{\rho})^2 \norm{\nb h}_{L^3(B_\rho)}
  \le N (\frac{r}{\rho})^2 [ \norm{\nb u}_{L^3(B_\rho)}
  + \norm{\nb v}_{L^3(B_\rho)} ] \\
\label{harmonic}
  &\le N (\frac{r}{\rho})^2 \norm{\nb u}_{L^3(B_\rho)}.
\end{align}
Combine this estimate with the above estimate for
$\norm{\nb v}_{L^3}$, divide by $r$,
and integrate in time to get
\begin{equation}
\label{tint}
  G_1(r) \le N \frac{\rho}{r} W(\rho)
  + N \frac{r}{\rho} G_1(\rho)
  + N \int_{-r^2}^0 |(\nb h)_r| dt.
\end{equation}
Since $h$ (and hence $\nb h$) is harmonic in $B_{3\rho/4}$,
$(\nb h)_r = (\nb h)_{\rho/2}$, and so
\[
  |(\nabla h)_r| = |(\nabla h)_{\rho/2}|
  \le \frac N\rho \norm{\nabla h}_{L^3(B_{\rho/2})}.
\]
Thus
\[
  g(u;r) := N \int_{-r^2} ^0 |(\nabla h)_r| dt
  \le \frac{N}{\rho}
  \int_{-r^2}^0  \norm{\nabla u}_{L^3(B_\rho)} dt .
\]
Since $\nb u \in L^{3,1}(Q_\rho)$,
we have $g(u,r) \to 0$ as $r \to 0$,
and so~\eqref{tint} yields~\eqref{wv-101}.
\end{pf}

\begin{rem} \label{RK3.7} By similar argument,
if $w = \curl u \in L^{p,q}_{loc}$ near $z$,
then so is $\nb u$, since
$\|\nabla u\|_{L^{p,q}(Q_r)} \leq
N \|w\|_{L^{p,q}(Q_{\rho})} + N \|u\|_{L^{p,q}(Q_{\rho})}$
if $0<r<\rho\le 2r$.
\end{rem}

\MAINPF{(iii)}

It suffices to prove the borderline cases $3/p+2/q=3$ and
$1 < p \le 3$. The other cases follow by H\"older's inequality.
If $p < 3$, we use the estimate~\eqref{wv-100},
and if $p=3$, we use the refined estimate~\eqref{wv-101}.
Choose $\th \in (0,1/4)$ so that if
$p < 3$, then $N\th^{3/p-1} < 1/2$, where $N$ is the constant
in~\eqref{wv-100}, and if $p = 3$, $N \th < 1/2$, where
$N$ is the constant from~\eqref{wv-101}.
Replace $r,\rho$ by $\th r$ and $r$, respectively.
Note that $G_1(r)$ is finite by Remark~\ref{RK3.7}.
The estimate \eqref{wv-100} ($p < 3$) or
\eqref{wv-101} ($p=3$) then implies
\[
  G_1(\th r) \leq \frac{N}{\th}W(r) + \frac{1}{2}G_1(r)
  + \left \{ \begin{array}{c}
  0 \;\; \mbox{ if } \;\; p < 3 \\
  N g(u; \theta r) \;\; \mbox{ if } \;\; p=3
  \end{array} \right..
\]
Choose $r_0$ so that $\sup _{r<r_0}W(r)<\frac{\th \e}{8N}$,
and (if $p=3$) $g(u;r_0) < \frac{\e}{8N}$,
where $\e$ is the constant in Theorem \ref{TH1} (ii).
Then for $r \le r_0$, we have
\[
  G_1(\theta r) \leq \frac{1}{2} G_1(r) + \frac{\e}{4}.
\]
Iterating this estimate, we obtain, for all $r \le r_0$,
\begin{equation}
G_1(\th^k r) \leq \frac{1}{2^k} G_1(r) +\frac{\e}2.
\label{vorticity-1000}
\end{equation}
Choose an integer $k_0 \ge 3 + \sup_{\th r_0 < r<r_0} \log_2 \frac{
G_1(r)}{\e}$.  Then $G_1(r) \le \e$ for $r<\th^{k_0}r_0$.  The
regularity of $u$ at $z=(0,0)$ now
follows from Theorem \ref{TH1} (ii).
\qed

\bigskip

In the next lemma, we show that vorticity and the
gradient of velocity can be controlled by the
gradient of vorticity in scaled norms, which is the key for
Theorem \ref{TH1} (iv).
\begin{lemm} \label{lemma3-7}
Suppose $0<2r\le\rho$ and $Q_{\rho}\subset \Om \times I$.
Suppose $1 \leq q \leq 2$, and $p, p^\sharp$ satisfy
\eqref{pq-100} and \eqref{grad-vorticity100}.
If $\nabla w\in L^{p^\sharp,q}_{x,t}(Q_{\rho})$, then
\begin{equation} \label{vorticity-grad-200}
  W(r) \leq N \bke{\frac{\rho}{r}}W_1(\rho)
  +N \bke{\frac{r}{\rho}}^{\frac{3}{p}-1}W(\rho).
\end{equation}
Furthermore, if $1\le q < 2$, we have
\begin{equation} \label{grad-vorticitycurl}
  G_1(r) \leq N \bke{\frac{\rho}{r}} \tilde{W}_1(\rho)
  + N \bke{\frac{r}{\rho}}^{3/p} G_1(\rho)
  + g(u;r)
\end{equation}
with $g(u;r) \to 0$ as $r \to 0$.
\end{lemm}
\begin{pf}
Statement~\eqref{vorticity-grad-200} follows from
Sobolev imbedding:
\begin{align*}
\norm{w}_{L^p(B_r)} &\le  \norm{(w)_{\rho}}_{L^p(B_r)} + \norm{w-
(w)_{\rho}}_{L^p(B_{\rho})}
\\
&\leq N\bke{\frac r\rho}^{\frac 3p} \norm{w}_{L^p(B_{\rho})} +
N\norm{\nabla w}_{L^{p^\sharp}(B_\rho)}.
\end{align*}
Taking the $L^q$ norm in time and dividing both sides by $r$,
we get~\eqref{vorticity-grad-200}.

The proof of~\eqref{grad-vorticitycurl}
is similar to that of the second part of Lemma~\ref{TH3-4}.
Choose a standard cut off function $\phi$ supported in
$B_{\rho}$ such that $\phi=1$ in $B_{3\rho/4}$. Define 2-tensors
\[
  V(x,t) := \int_{\R^3} \nabla_x \frac{4\pi}{|x-y|}
  \bke{\curl w(y,t)}\phi(y)dy,
  \quad\quad H := \nabla u - V.
\]
Note that since $\nabla \cdot u = 0$, $\Delta u = \curl w$, and so
$\Delta_x H(x,t) = 0$ in $B_{3\rho/4}$.
Potential estimates give, if $p ^\sharp >1$, (i.e.~$q<2$),
\[
  \norm{ V }_{L^p} \le N \norm{\curl w}_{L^{p^\sharp}(B_\rho)}.
\]
The same estimate as in~\eqref{harmonic} gives
\[
  \norm{H - (H)_r}_{L^p(B_r)}
  \le N(\frac{r}{\rho})^{\frac 3{p^\sharp}} \norm{H}_{L^p(B_\rho)}
  \le N(\frac{r}{\rho})^{\frac 3{p^\sharp}}
  \bkt{\norm{\nb u}_{L^p(B_\rho)} + \norm{\curl w}_{L^{p^\sharp}(B_\rho)}}.
\]
Using $\frac 1r \norm{(H)_r}_{L^p(B_r)} = N r^{\frac 3p -1} |(H)_r| =
N r^{\frac 3p -1} |H|_{r=0}$ together with the last two estimates, we
find \eqref{grad-vorticitycurl}
with
\[
g(u;r) :=  N r^{\frac3p-1} (\int_{-r^2}^0 (|H|_{r=0})^{q}\ind
t)^{1/q}.
\]
Now arguing as in the proof of Lemma~\ref{TH3-4},
\[
  g(u;r) 
  = N r^{\frac3p-1}(\int_{-r^2}^0 |(H)_{\rho/2}|^q dt)^{1/q}
  \le  \frac {Nr^{\frac3p-1}}{\rho^{3/p}}
  (\int_{-r^2}^0  (\norm{\nabla u}_{L^p(B_\rho)})^q dt)^{1/q},
\]
and since $\nb u \in L^{p,q}(Q_{\rho})$,
we have $g(u,r) \to 0$ as $r \to 0$.
\end{pf}

\MAINPF{(iv)}

It suffices to prove the borderline cases $3/p^\sharp +2/q=4$.
The other cases follow by H\"older inequality.  We also assume $z=(0,0)$.

We first consider $1 < q \leq 2$.
Let $\e$ be the constant in Theorem
\ref{TH1} (iii) and we suppose $W_1(r) \le \e/4$ for any $r<r_0$.
Our assertion follows a procedure similar to the proof in
Theorem \ref{TH1} (iii).
Since $\frac 3p-1 >0$ in \eqref{vorticity-grad-200}, we replace $r$ and
$\rho$ by $\theta r$ and $r$, respectively, after choosing
$\theta\in (0,1/2)$ appropriately, and then iterate
\eqref{vorticity-grad-200}. This procedure leads to the
conclusion that $W(\th^k r) < \e/2$ for $r<r_0$ and $k$
sufficiently large, and so $W(r)<\e$ for $r<\th^k r_0$.
Theorem \ref{TH1} (iii) then implies the regularity.

For $1 \le q<2$, we can use~\eqref{grad-vorticitycurl} instead.
Arguing just as in the second part of the proof of Theorem~\ref{TH1}(iii),
we conclude that $G_1(r)$ can be made small enough to apply
Theorem~\ref{TH1}(ii), provided $\tilde{W}_1(r) \leq W_1(r)$
can be made arbitrarily small.
\qed


\section*{Acknowledgments}

{The research of Gustafson and Tsai is partly supported by NSERC
grants.  Kang thanks the University of British Columbia and the
Pacific Institute of Mathematical Sciences for their hospitality
during his visit in February, 2006.}


\end{document}